\magnification=1200

\font\smallcaps=cmcsc10

\def\qed{${\vcenter{\vbox{\hrule height .4pt
           \hbox{\vrule width .4pt height 4pt
            \kern 4pt \vrule width .4pt}
             \hrule height .4pt}}}$}

\def\mqed{{\vcenter{\vbox{\hrule height .4pt
           \hbox{\vrule width .4pt height 4pt
            \kern 4pt \vrule width .4pt}
             \hrule height .4pt}}}}


\def\CC{\rlap {\raise 0.4ex \hbox{$\scriptscriptstyle |$}}
  \hskip -0.1em C}
\def\FF{\hbox to 8.33887pt{\rm I\hskip-1.8pt F}}
\def\NN{\hbox to 9.3111pt{\rm I\hskip-1.8pt N}}
\def\PP{\hbox to 8.61664pt{\rm I\hskip-1.8pt P}}
\def\QQ{\rlap {\raise 0.4ex \hbox{$\scriptscriptstyle |$}}
  {\hskip -0.1em Q}}
\def\RR{\hbox to 9.1722pt{\rm I\hskip-1.8pt R}}
\def\ZZ{\hbox to 8.2222pt{\rm Z\hskip-4pt \rm Z}} 

\centerline{Counting Points on Curves over Families in Polynomial Time}
\bigskip
\centerline{ Jonathan Pila }
\bigskip
\medskip

\noindent
The purpose of this note is to apply the result of [9] to obtain an
algorithm for counting points on curves over finite fields,
for curves that belong to families.  By a family of curves we mean
curves that are specified by equations of the same form, in 
an ambient projective space of the same dimension, and have the same
degree and genus. We will be more precise below. The algorithm
will be polynomial time in the size of the finite field, with the degree
of the polynomial depending on the family.

The main result of [9] was the following generalization
to Abelian varieties over finite fields of
the algorithm of Schoof [10] for elliptic curves over finite fields. 

\medbreak
{\bf Theorem A.\ }{\it Let $A$ be an Abelian variety over a finite field
${\bf F}_q$,
given explicitly as a projective variety with an explicit addition law.
Then one can compute the characteristic polynomial of the Frobenius
endomorphism of $A$ in time bounded by $B_1(\log q)^{B_2}$, 
where $B_1,B_2$ 
depend only on the embedding space dimension of $A$, the
number of equations defining $A$ and the addition law, and their degrees.\ }
\qed
\medbreak

Let us remark at this point that the above is not a generalization of
Schoof's result in the strict sense: the elliptic curve case of the above
does not reduce to Schoof's algorithm. In particular, we require for input
a more comprehensive description of the group law (see [8, section 3]). 
The group law on an
elliptic curve is traditionally specified by a set of rational functions
defined on an open subset of $E\times E$; another set of rational
functions gives the group law on the complement of that open set in
$E\times E$. Such a specification suffices
for Schoof. We require that several such sets of rational functions
be given, each determining the
group law on an open subset of $E\times E$, and such that these open
subsets cover $E\times E$. See Silverman [12, Remark 3.6.1] 
for a discussion of how
to obtain further charts of the addition law in the elliptic curve case.

In order to apply our algorithm to count points on curves over finite
fields, we must first construct the Jacobian variety of the curve,
and the group law upon it in the above sense.
In the applications in [9], this was done following the construction of
Chow [3]. Applying Chow's construction to a curve $C$ defined over {\bf Q},
one obtains in explicit form the Jacobian $J$.
One can then obtain the Jacobian of the reduction of $C$ modulo
primes $p$ (for almost all $p$) by reducing $J$.
These reductions are all defined by equations of the same form.

These applications were thus of the type considered here: the family of
curves being the reductions mod $p$ of a fixed curve over {\bf Q}. They
depended on the compatibility of Chow's construction with reduction.

However, Chow's construction is also compatible with other
specializations; the following is a theorem of Igusa [7], who also
analysed the situation in which the specialization of the curve
acquires some limited singularities.
A still more general statement of the universality of the Jacobian
construction is due to Grothendieck [6]. A convenient reference is
Milne [8, section 8]. The statement below is sufficient for our purposes.

\medbreak
{\bf Theorem B.\ }{\it Let $C$ and $C'$ be two irreducible curves in
${\bf P}^n$ such that $C'$ is a specialization of $C$ over $K$.
Let $J$ and $J'$ be the completed generalized Jacobian varieties of
$C$ and $C'$ respectively.
Suppose that $C$ and $C'$ have the same arithmetic genus, and that
the same reference integer is used in the constructions of the Jacobian
varieties. Then $J$ and $J'$ have the same ambient space.
If $C$ is nonsingular, and if $C'$ is either
nonsingular or has one and only one ordinary double point, then $J'$ is the
unique specialization of $J$ over the specialization $C\rightarrow C'$ with
reference to $K$.\ \qed}
\medbreak

To state a result on the form of the equations defining the Jacobian, we
begin with a definition.

\medbreak 
{\bf Definition.\ } A {\it curve presentation\ } is a finite collection
of forms (in the variables $X_i$) in
$${\bf Z}[a_1,\ldots, a_m][X_0,\ldots, X_n]$$
where the $a_j$ are indeterminates,
with the property that if the $a_j$ are considered as independent
transcendental elements over {\bf Q}, the projective algebraic set determined
by the forms is an absolutely irreducible nonsingular curve. We refer to
the degree and genus of this curve as the degree and genus of the
presentation.
\medbreak

{\bf Remark.\ }
The restriction that the $a_j$ be independent indeterminates is not
necessary. One can allow them to be connected by some arbitrary
algebraic relations. This gives a family of curves corresponding to the
points of some algebraic set, $V$, given as the zeros of an ideal $\Lambda$. 
The stated definition takes $V$ to be affine $m$-space. 

If the ideal $\Lambda$ is prime
then we require that the curve determined by the
forms over the field of fractions of {\bf Z}$[a_j]/\Lambda$ (that is, the
generic point)
be an absolutely irreducible non-singular
curve. If $V$ is reducible, we require the above at all constituent primes.
For simplicity we will assume that our presentations are irreducible;
reducible presentations can be dealt with by first finding (in a finite
amount of time [11]) the constituent irreducible presentations.

\medbreak

Associated with a curve presentation {\bf C\ }is a finite collection 
$\Delta({\bf C})$ of polynomials $\Delta_i$,
$$\Delta_i \in {\bf Z}[a_1,\ldots, a_m],$$
with the following property: if the $a_j$ are specialized to values in an
arbitrary field $K$, the projective algebraic set determined by the forms
of the presentation will be an irreducible nonsingular curve,
having the degree and
genus of the presentation, except possibly when the polynomials $\Delta_i$
all vanish.
We call $\Delta({\bf C})$ the {\it discriminant\ } of the presentation.
\medbreak
To build the identity element of the Jacobian variety we need some
rational divisor on the curve presentation.
The identity element will be a suitable multiple of this divisor. 
This should be
specified as being a point (with coordinates in ${\bf Z}[a_1,\ldots, a_m]$)
in some fixed symmetric product of the curve presentation. For simplicity
we will assume that it is given as the intersection of the curve with
a hyperplane.
\medbreak
{\bf Definition.} A {\it rational divisor presentation\ }{\bf D} 
associated to
a  curve presentation {\bf C\ }is a linear form in 
${\bf Z}[a_1,\ldots ,a_m][X_0,\ldots ,X_n]$ with the property that
when the $a_i$ as considered as independent transcendental elements
over {\bf Q}, the hyperplane determined by the form properly intersects
the curve.
\medbreak
Associated with a curve presentation and a rational divisor
presentation is a discriminant $\Delta ({\bf C, D})$ with the property
that when the $a_i$ are specialized to values in an arbitrary field,
the specialization of the hyperplane is a hyperplane (that is, not
all the coefficients vanish), and this hyperplane properly intersects
the specialization of the curve, except possibly when the polynomials
comprising the discriminant all vanish.

\medbreak
{\bf Definition.\ }An {\it Abelian variety presentation\ }{\bf A\ }
consists of the 
following. Let $a_j$ be indeterminates.

1. A finite collection of forms $F_i(X)$ in 
${\bf Z}[a_1,\ldots ,a_m][X_0,\ldots ,X_n];$

2. A finite collection of $n+1$-tuples of polynomials 
$$G^{(r)}(X,Y),\;\; r=1,\ldots ,R$$
$$G^{(r)}(X,Y)=(G^{(r)}_0(X,Y),\ldots ,G^{(r)}_n(X,Y)),$$
where $X=(X_0,\ldots ,X_n), Y=(Y_0,\ldots ,Y_n)$ and the
$G^{(r)}_i$ are homogeneous of the same degree in each system of variables.

3. An $n$-tuple $E$ of elements of ${\bf Z}[a_1,\ldots ,a_m]$, not all
zero.

\noindent
These objects should further have the property that when the $a_j$ are 
considered as independent transcendental elements over {\bf Q\ }, the
forms $F_i$ determine a nonsingular variety $A$, 
the $n$-tuple $E$ point of this
variety, and the collection of $n+1$-tuples $G^{(r)}$ a group law 
$A\times A \rightarrow A$ with the point $E$ as identity element.
By this last condition we mean that each $n+1$-tuple should define
the group law on an open subset of $A\times A$, and that these open subsets
should together cover $A\times A$.

Here again we will allow the $a_j$ to be connected by some
algebraic relations. We will insist that the corresponding
ideal of relations be prime.
\medbreak
Associated with an Abelian variety presentation is a discriminant
$\Delta({\bf A})$ comprising a finite number of polynomials
in ${\bf Z}[a_1,\ldots ,a_m]$. It has the
property that if the $a_j$ are specialized to values in an arbitrary field,
the corresponding algebraic set is an abelian variety, with
$E$ a point of this variety, with group law determined by the forms
$G^{(r)}$ in the above sense, and with $E$ as the identity element,
except possibly when the constituent polynomials
of the discriminant all vanish.
\medbreak
Given two systems of polynomials $\Delta_1$ and $\Delta_2$, the product
$\Delta_1\Delta_2$ will denote the system of generators for the
product of the ideals.
We can now state an equational version of 
Chow's construction and Igusa's theorem.
\medbreak
{\bf Theorem C.\ }{\it Let {\bf C} be a curve presentation, and
{\bf D} an associated rational divisor presentation. Then we can
construct an Abelian variety presentation {\bf J} 
and a finite non-empty collection of polynomials $\Delta^*$ with the
following property. Let the $a_j$ be specialized to values
in an arbitrary field in such a way that 
the polynomials
$\Delta^* \Delta ({\bf C}) \Delta ({\bf C},{\bf D})$
do not all vanish. Then the Abelian variety
presentation determines the Jacobian variety of the curve.}
\medbreak
We call {\bf J\ }the {\it Jacobian presentation\ } associated to
the curve presentation {\bf C}, and the rational divisor presentation
{\bf D}. 
It is completely determined by
the curve presentation, the rational divisor presentation, and the
choice of an auxilliary positive integer in the construction.
The discriminant $\Delta({\bf J})$ thus divides 
$\Delta^*\,\Delta({\bf C})\,\Delta({\bf C},{\bf D})$.

\medbreak
{\bf Proof.\ }
Chow's construction
obtains equations for the variety underlying the Jacobian, and for
the graph of the group law, in universal form
for the presentation that specialize whenever the curve and rational divisor
specialize appropriately. We desire a rational 
parametrization of this graph.
Initially, we construct such a parametrization over {\bf Z}$[a_j]$.
The forms $G^{(r)}_i$ we construct will 
give a complete description of the group law on $J\times J$ over
${\bf Q}(a_1,\ldots ,a_m)$ when the $a_j$ are taken to be
independent transcendental elements; it follows that the forms $G^{(r)}_i$
have no common zeros on $J\times J$. By the nullstellensatz, the ideal
generated by the $G^{(r)}_i$ and the $F_i$ over {\bf Z\ }$[a_j]$
contains a non-zero element of {\bf Z\ }$[a_j]$;
moreover, by the elimination theorem, there are a finite number
of elements of {\bf Z\ }$[a_j]$ whose vanishing is a necessary
and sufficient
condition for the existence of a common zero (in the $X_i$) of
$G^{(r)}_i$ and the $F_i$.
This non-empty collection of non-zero polynomials gives the
additional discriminant $\Delta^*$ of the theorem.\ \qed
\medbreak

We thus have a presentation of the Jacobian that presents
the Jacobian in a form suitable for the application of theorem A,
except on the zero set of $\Delta^*$, a lower dimensional set.  
On the zero set of $\Delta^*$ (but outside that of $\Delta$) we have a
presentation of the variety, and the graph of the group law.
We can determine ([11]) the associated primes of  $\Delta^*$.
On each, we can consider the Abelian variety
over the corresponding field of fractions, and construct a
parametrization of the graph of the group law.
Repeating the argument of the above proof, the forms we construct will
continue to give a complete description of the group law when
specialized outside a set of lower dimension, determined by
some resultant system. In this way we can construct a finite number
of Abelian variety presentations and systems of discriminants
such that, for any specialization of the curve and divisor, outside their
discriminants, one (at least) of the Abelian variety presentations 
will present the Jacobian of the curve. (Since
the irreducible components of a given algebraic set will intersect,
a given specialization of the $a_j$ might correspond to more than one
of the Abelian variety presentations.)

We thus conclude the following result.

\medbreak
{\bf Theorem D.\ }{\it Let {\bf C} be a curve presentation, and {\bf D}
an associated rational divisor presentation. 
Let $\Delta$ be $\Delta({\bf C})\,\Delta({\bf C},{\bf D})$.
Then we can construct
a finite number of Abelian variety presentations ${\bf J}_1,\ldots,
{\bf J}_r$, corresponding to prime ideals $\Lambda_i$ in ${\bf Z}[a_j]$,
and with discriminants $\Delta_i$, with the following property.
Let the $a_j$ be specialized to values in any field in such a way that the
polynomials comprising $\Delta$ do not all vanish. Then one (at least) of
the Abelian variety presentations ${\bf J}_i$
will have the property that the
$a_j$ belong to the corresponding set $V(\Lambda_i)$, and the polynomials
comprising the corresponding
discriminant $\Delta_i$ do not all vanish.
Further, any of the presentations with the above property will
determine the Jacobian variety of the curve.\ \qed
}
\medbreak
{\bf Remark.\ }Let us raise the question of whether the above
stratification procedure is necessary: perhaps the kind of data we desire
describing the group law can be obtained universally, but
we do not see how to do this.

More generally, suppose that $V$ and $W$ each represent smooth, complete
varieties, parametrized by points of another variety (outside
some discriminant), and determined by equations in a universal way. 
Suppose further that we have a universal presentation of the graph
of a morphism $V\rightarrow W$. Can we get a universal system of rational
functions determining this morphism? 

If the discriminant of the original presentation of the curve is
of codimension 1, one might expect that the discriminant $\Delta^*$
already coincides with $\Delta$.
\medbreak
Let us now eliminate the divisor presentation from the data.
Given a curve presentation, we can find some rational divisor
presentation that generically intersects the curve
properly. 
We now proceed by a similar process of stratification to
determine divisors on the lower dimensional presentations.
We thus get a version of theorem D with no rational divisor in the
hypothesis.
\medbreak
{\bf Theorem E.\ }{\it
Let {\bf C} be a curve presentation with discriminant
$\Delta$.
Then we can construct
a finite number of Abelian variety presentations ${\bf J}_1,\ldots,
{\bf J}_s$, corresponding to prime ideals $\Lambda_i$ in ${\bf Z}[a_j]$,
and with discriminants $\Delta_i$, with the following property.
Let the $a_j$ be specialized to values in any field in such a way that the
polynomials comprising $\Delta$ do not all vanish. Then one (at least) of
the Abelian variety presentations ${\bf J}_i$
will have the property that the
$a_j$ belong to the set 
corresponding to $\Lambda_i$, but not to the set
corresponding to $\Delta_i$. 
Further, any of the presentations with the above property will
determine the Jacobian variety of the curve.\ \qed
}
\medbreak
Combining with theorem A yields the following algorithm.

\medbreak
{\bf Theorem F.\ }{\it Let {\bf C\ }be a curve presentation
of degree $d$ and genus $g$.  
There exist positive integers $B_1,\; B_2$, 
and a deterministic algorithm,
depending only
on the presentation,  with the following
property. Let the $a_j$ be specialized to values in a finite
field ${\bf F}_q$ such that the 
polynomials comprising $\Delta ({\bf C})$ do not all vanish.
Then the algorithm computes the Zeta function of the curve over
the finite field ${\bf F}_q$, and hence in particular
the number of rational points on the curve over ${\bf F}_q$ 
in time bounded by
$B_1\; (\log q)^{B_2}.$\ \qed\ }
\medbreak

As an example of the above we will consider the hyperelliptic
families 
$$y^2\; =\; f(x)\; =\; a_0 x^d+a_1x^{d-1}+\cdots +a_d$$
where $d\ge 3$ is an integer. This is not a presentation according to
our definition, since the equation is not homogeneous, and the
plane curve is singular (if $d\ge 5$). However, a (completed) 
projective
embedding is easily given, in universal form (see Silverman [11, exercise
2.14]). The curve is nonsingular of degree $d$ and genus $g$, where
$d-2\le 2g < d $ for a given specialization of the $a_j$
provided that the discriminant $\Delta (d)$ of the polynomial $f(x)$
is non-vanishing. 
In this case a rational divisor presentation is easily found
with discriminant 1. Hence we obtain the following result.

\medskip
{\bf Theorem G.\ }{\it Let $d\ge 3$ be an integer. There exist positive
integers $B_1, B_2$,
and a deterministic algorithm, depending only on $d$,
with the following property. Let  $f(x)$ have
coefficients $a_j$ in a finite field
${\bf F}_q$, such that the discriminant $\Delta (d)$
does not vanish.
Then the
algorithm computes the number of points on the corresponding
hyperelliptic curve over ${\bf F}_q$ in time bounded by 
$B_1\; (\log q)^{B_2}.$\ \qed\ }
\medbreak

A random polynomial time algorithm in the genus 2 case has been given by
Adleman and Huang [1] as part of their random polynomial time primality test.

Let us also remark that the above, and all the algorithms presented here
are of purely theoretical interest (if any) at the present time: not
only because the algorithm of theorem A is impractical, but because
the construction of the Jacobian in the desired form is impractical.
As remarked earlier, this is generally avoided even in the elliptic curve
case.
\medbreak
As a further application we can give a uniform version of theorem C of [9].
For this, we consider initially a 
family of projective plane curves, generically absolutely
irreducible, of a certain genus and degree. 
We allow a finite number of singularities. 
All these properties will be preserved by specializations outside
an appropriate discriminant. 
Such a family is specified by a single form $H(X,Y,Z)$ in
{\bf Z}$[a_1,\ldots ,a_m][X, Y, Z]$,
with some other relations on the $a_j$ given by
a (prime) ideal $\Lambda$ in ${\bf Z}[a_j]$.

\medbreak
{\bf Theorem H.\ }{\it Let {\bf H} be a family of plane curves,
with discriminant $\Delta$.
There is a deterministic polynomial time algorithm operating as follows.
Taking as input a specialization of the $a_j$ in an arbitrary
finite field ${\bf F}_q$
such that the polynomials in $\Lambda$ all vanish,
but not all those in $\Delta$, the algorithm counts the number of
points on the curve $H(X,Y,Z)=0$ in ${\bf P}^2({\bf F}_q)$
}
\medbreak
{\bf Proof.\ } We show how to obtain a finite number of curve
presentations from the planar family. 
We begin by resolving the generic curve of the family.
This gives a smooth projective curve presentation, outside
of a set of lower dimension determined by
an appropriate discriminant. 
After finding the constituent primes of the zero set of this
discriminant, one resolves the corresponding curves.
In this way, one produces a finite number of curve presentations
that suffice to resolve every curve in the family.
We can now apply
theorem F. The number of rational points will be the same except for
a correction due to the singular points.
The correction is easily made in polynomial time, as it entails
counting points in a zero dimensional algebraic set with a fixed
description (see [9]).
\ \qed
\medbreak
One can state a similar theorem for families of affine plane curves.
The extra accounting for the points at infinity is easily accomplished
in polynomial time.
\medbreak
We now make some concluding remarks.
Apart from using simpler input data, Schoof is able to reduce
the computations needed to calculations involving univariate
polynomials (essentially by operating on the Kummer variety).
The algorithm of theorem A does ideal-theoretic calculations
with polynomials in many variables.
In this sense our result might be paraphrased as the
assertion that Schoof's idea of computing using the
$\ell$-adic representations is strong enough to succeed without
further assistance from simplifications in the ideal-theoretic
computations.

Nevertheless one would like to investigate such simplifications in the
higher dimensional cases, especially the genus or dimension 2 case
that is the next simplest after elliptic curves. Explicit equations
for Jacobian varieties of hyperelliptic curves of genus two,
together with group laws (seemingly not in the comprehensive form we require)
are now available, due to Grant [5] (degree 5), and
Cassels and Flynn [2], [4] (degree 6). These provide a starting point.

Returning to theoretical questions, it is natural to pursue an algorithm of
the type of theorem E applicable to a family containing all curves
of a given genus. 

In order to avoid subtleties 
connected with the moduli space, it might be desirable
to work instead with a family that represents {\it most\ }isomorphism
classes of curves, each a {\it roughly\ }equal number of times.
Thus we could (say) leave out hyperelliptic curves ($g\ge 3$), and also
allow curves with automorphisms to be under-represented. The object
would be a family of curves imitating the moduli space in a
statistical sense: with respect to the distribution of
coefficients of Zeta functions over finite fields. A similar question
arises for Abelian varieties of a given dimension.

\bigskip
\noindent
{\bf Note (April 2005).\/} This paper was written in (March) 1991,
and I have posted it ``as was''. At a later point I intend to revise
it and update references. My thanks to Claus Diem for his suggestions.
At the time the paper was written, I was
affiliated with the Department of Mathematics of Columbia University
(New York).

\bigskip
\bigskip
\bigskip
\noindent
{\bf References.}
\medskip
\frenchspacing

\item{[1]} {\smallcaps Adleman,\thinspace L.\thinspace M. 
and Huang,\thinspace M.--D.,} 
Recognizing Primes in Random Polynomial Time, Preprint,  1988.

\item{[2]} {\smallcaps Cassels,\thinspace J.\thinspace W.\thinspace S.,}
Arithmetic of curves of Genus 2, {\it Number Theory and Applications,}
R. A. Mollin, editor, Kluwer 1989.

\item{[3]} {\smallcaps Chow,\thinspace W.--L.,}  
The Jacobian variety of an algebraic curve, 
{\it American. J. Math.} {\bf 76} (1954), 453--476.

\item{[4]} {\smallcaps Flynn, \thinspace E.\thinspace V.,}
The jacobian and Formal Group of a Curve of Genus 2 over an Arbitrary
Ground Field, 
{\it Math. Proc. Camb Phil. Soc.} {\bf 107} (1990), 425--441.

\item{[5]} {\smallcaps Grant,\thinspace D.,}
Formal Groups in Genus Two,
{\it J. Reine Angew. Math.,} {\bf 411} (1990) 96--121.

\item{[6]} {\smallcaps Grothendieck,\thinspace A.,\ }
Technique de descente et th\'eor{\accent18 e}ms d'existence\ en g\'eom\'etrie
alg\'ebrique, IV,V,
{\it S\'eminaire Bourbaki,} \'Expos\'es 221, 232. (1960-1962).

\item{[7]} {\smallcaps Igusa,\thinspace J--I.,}
Fibre Systems of Jacobian Varieties,
{\it American J. Math.,} {\bf 78} (1956), 171--199.

\item{[8]} {\smallcaps Milne,\thinspace J.\thinspace  S.,}
Jacobian Varieties, in
{\it Arithmetic geometry,} Cornell and Silverman editors,
Springer Verlag, New York 1986.

\item{[9]} {\smallcaps Pila,\thinspace J.,}
Frobenius Maps of Abelian Varieties and Finding Roots of Unity in 
Finite Fields,
{\it Math. Comp.} {\bf 55} (1990), 745--763. 

\item{[10]} {\smallcaps Schoof,\thinspace R.\thinspace J.,} 
Elliptic curves over finite fields and the computation of square 
roots mod
$p$, {\it Math. of Comp.} {\bf 44} (1985), 483--494.

\item{[11]} {\smallcaps Seidenberg,\thinspace A.,}
Constructions in Algebra,
{\it Trans. Amer. Math. Soc.,}{\bf 197} (1974),273--313.

\item{[12]} {\smallcaps Silverman, \thinspace J.\thinspace H.,} 
{\it The Arithmetic of 
Elliptic  Curves,} Springer Verlag, New York 1986.

\nonfrenchspacing

\bigskip
\bigskip
\leftline{
Department of Mathematics and Statistics, McGill University, 
Montreal H3A 2K6 Canada}
\bigskip
\noindent
pila@math.mcgill.ca

\bye